\documentclass[12pt]{article}
\usepackage{amsmath,amssymb,amsfonts}

\newcommand{\D}{\displaystyle}
\newcommand{\R}{\mathbb{R}}
\newcommand{\N}{\mathbb{N}}
\newcommand{\lo}{\longrightarrow}
\newcommand{\ov}{\overline}
\newcommand{\para}{\paragraph}
\newcommand{\var}{\varepsilon}
\begin{document}
\begin{center}
{\bf \large Palais-Smale Condition, Index Pairs and Critical Point
Theory }\\[1cm] M.R. Razvan
\\ \noindent {\it \small Institute for Studies in Theoretical Physics and
Mathematics\\ \small P.O.Box: $19395-5746$, Tehran, IRAN}\\
e-mail: razvan@karun.ipm.ac.ir \\ Fax: 009821-2290648
\end{center}

\begin{abstract}
This paper is concerned with index pairs in the sense of Conley
index theory for flows relative to pseudo-gradient vector fields
for $C^1$-functions satisfying Palais-Smale condition. We prove a
deformation theorem for such index pairs to obtain a
Lusternik-Schnirelmann type result in Conley index theory.
\end{abstract}

\noindent {\bf Keywords:} Critical point theory, index pair,
Palais-Smale condition, pseudo-gradient vector field,  relative
Lusternik-Schnirelmann category. \\ {\bf Subject Classification:}
37B30, 58E05

\section{Introduction}
Conley's homotopy index has proved to be a useful tool in critical
point theory. In most of the applications of critical point
theory, we deal with a $C^1$-function on a complete Finsler
manifold satisfying Palais-smale condition and we desire to
estimate the number of critical points of such a function. There
are two different viewpoints: the number of analytically distinct
critical points which is discussed by Morse theory and the number
of geometrically distinct critical points which is investigated by
Lusternik-Schnirelmann theory \cite{LS}. The applications of
Conley index in Morse theory has been studied in details
\cite{B1,B2,BD,BG,CZ2}. This paper is concerned with the
applications in Lusternik-Schnirelmann theory. (See also
\cite{Po,R1,R2,R3}.) We consider the flow relative to a
pseudo-gradient vector field for a $C^1$-function on a complete
Finsler manifold satisfying Palais-Smale condition. Then we prove
the deformation results for
 index pairs in the sense of Conley index theory for such a flow. Using a
modification of the relative Lusternik-Schnirelmann category, we
obtain a Lusternik-Schnirelmann type result in Conley index
theory.

We need some basic results from Conley index theory which are
presented in Section 2. Here we have extended the concept of
regular index pair \cite{S} for the noncompact case. Then in
Section 3, we define the concept of relative
Lusternik-Schnirelmann category with a little modification so
that it can be applied to index pairs and yield sharper results
in critical point theory. Finally in Section 4, we prove our
result concerning the existence of critical points for a
$C^1$-function on a complete Finsler manifold.

\section{Index Pair}
Let $\varphi^t$ be a continuous flow on a metric space $X$. An
{\em isolated invariant set} is a subset $I\subset X$ which is the
maximal invariant subset in a closed neighborhood of itself. such
a neighborhood is called an {\em isolating neighborhood}. In order
to define the concept of index pair, we follow \cite{B1} and
\cite{RS}.

Given a closed pair $(N,L)$ in $X$, we define the induced
semi-flow on $N/L$ by $$\varphi_{\sharp}^t: N/L \lo N/L , \ \
\varphi_{\sharp}^t (x)=\left\{
\begin{array}{ll}
\varphi^t (x) & \text{if} \ \varphi^{[0,t]} (x) \subset N-L\\
\left[ L \right] & \text{otherwise.}
\end{array} \right. $$
In \cite{R2,RS} it is proved that $\varphi _{\sharp}:N/L \times \R
^+ \lo N/L $ is continuous if and only if
\begin{itemize}
\item[(i)] $L$ is positively invariant relative to $N$, i.e. $$x\in L,
t\geq 0, \varphi^{[0,t]} (x) \subset N \Rightarrow \varphi^{[0,t]}
(x) \subset L.$$
\item[(ii)] Every orbit which exits $N$ goes through $L$ first, i.e.
$$x\in N,\varphi^{[0,\infty)} (x) \not\subset N \Rightarrow
\exists_{t\geq 0} \ \text{with} \ \varphi^{[0,t]} (x) \subset N,
\varphi^t (x) \in L,$$ or equivalently if $x\in N-L$ then there is
a $t>0$ such that $\varphi^{[0,t]} (x) \subset N$.
\end{itemize}
\para{Definition.} An {\em index pair} for an isolated invariant set
$I\subset X$ is a closed pair $(N,L)$ in $X$ such that
$\overline{N-L}$ is an isolating neighborhood for $I$ and the
semi-flow $\varphi^t_{\sharp}$ induced by $\varphi^t$ is
continuous.

For every subset $A\subset X$ and $T\in \R^+$, we define
$$\begin{array}{l} G^T(A)=\{x\in A|\varphi^{[-T,T]} (x) \subset
A\}, \\ \Gamma^T (A)=\{x\in G^T(A)|\varphi^{[0,T]} (x) \cap
\partial A\neq \varnothing\}.
\end{array}$$
Now suppose that $G^T(A)\subset int(A)$ for a closed subset
$A\subset X$. Then $A$ is an isolating neighborhood for
$I:=\bigcap_{T>0} G^T(A)$. In \cite{B1}, Benci proved that
$(G^T(A), \Gamma^T(A))$ is an index pair for $I$. This index pair
has a special property, that is for every $x\in \Gamma ^T(A)$,
$\varphi ^{[0,3T]}(x) \not \subset G^T(A)$. To see this, suppose
that $\varphi ^{[0,3T]} \subset G^T(A)$, then $\varphi
^{[-T,4T]}(x)\subset A$, hence $\varphi ^{[T,2T]}(x) \subset
G^{2T}(A)\subset int(G^T(A))$. So $\varphi ^T (x)\not \in \Gamma
^T(A)$ and by (i), $x\not \in \Gamma ^T(A)$.

\para{Definition.} An index pair $(N,L)$ is called {\em regular} if the
exit time map $$\tau_+:N \lo [0,+\infty] , \ \ \tau_+(x)=\left\{
\begin{array}{ll}
sup\{t|\varphi^{[0,t]} (x) \subset N-L\} & \text{if} \ x\in N-L,\\
0 & \text{if}\ x\in L,
\end{array} \right. $$
is continuous. For every regular index pair $(N,L)$, we define the
induced semi-flow on $N$ by $$ \varphi ^t _{\natural} :N\times \R
^+ \lo N,\ \ \varphi^t _{\natural}(x)=\varphi ^{min\{t,\tau_+
(x)\}} (x)$$

The reader is referred to \cite{S} to see the details about
regular index pairs and the proof of the following useful
criterion.

\para{Proposition 2.1.} An index pair $(N,L)$ is regular provided
that $\varphi ^{[0,t]}(x)\not \subset \ov{N-L}$ for every $x\in L$
and $t>0$.

\para{Definition.} An index pair $(N,L)$ is called {\em weakly
regular} if for every $x\in L$, there exists $t\in \R ^+$ such
that $\varphi ^t (x)\not \in \ov{N-L}$.

We showed that if $G^T(A)\subset int(A)$ for a closed subset
$A\subset X$, then $(G^T(A), \Gamma^T(A))$ is a weakly regular
index pair. Now let $(N,L)$ be an index pair and $V:=\ov{N-L}$. We
define $I^+(V)=\{x\in V|\varphi^{[0,+\infty)} (x) \subset V\}$.
It is not hard to see that $I^+ (V)$ is closed subset of $V$.
Notice that if $(N,L)$ is a weakly regular index pair, then
$I^+(V)\bigcap L= \varnothing$. The following theorem asserts
that there is a Lyapunov function for $\varphi ^t _{\sharp}$ on
$N/L$ which separates $[L]$ and $I^+(V)$. If we consider the
natural projection $\pi :N\lo N/L$, then we obtain a Lyapunov
function on $N$ which separates $L$ and $I^+(V)$. (See \cite{S}
for a similar result in the compact case.)

\para{Theorem 2.2.} Let $(N,L)$ be a weakly regular index pair and
$V=\ov{N-L}$. Then there exists a continuous function $g:N/L\lo
[0,1]$ such that
\begin{itemize}
\item[(i)] $g^{-1} (0)=[L]$ and $g^{-1}(1)=I^+ (V)$.
\item[(ii)] If $0<g(x)<1$ and $t\in \R^+$, then
$g(\varphi^t_{\sharp}(x))<g(x)$.
\end{itemize}
{\it Proof.} Let $\rho :N/L\lo [0,1]$ be a continuous function
with $\rho^{-1}(0)=[L]$ and $\rho^{-1}(1)=I^+(V)$. We define
$f:N/L\lo [0,1]$ by $f(x)=\D{\sup_{t\geq 0}}\ \rho
(\varphi^t_{\sharp} (x))$. Then $f^{-1} (0)=[L],\ f^{-1} (1)=I^+
(V)$ and $f(\varphi^t_{\sharp} (x))\leq f(x)$ for every $x\in
N/L$. Moreover it is not hand to check that $f$ is continuous in
$I^+(V)$. Now suppose that $0<f(x)<1$ for some $x\in N-L$. Then
there is $t_0\in \R^+$ such that $\varphi^{t_0} (x)\in L$. Since
$(N,L)$ is a weakly regular index pair there is a $t_1\in\R^+$
such that $\varphi^{t_0+t_1} \not\in V$. Thus there is a
neighborhood $A$ of $x$ such that $\varphi^{t_0+t_1} (A) \cap V
=\varnothing$ which means that $\varphi_{\sharp}^{t_0+t_1}
(A)=[L]$. Therefore $\D{\sup_{t\geq 0}}\ \rho(\varphi^t_{\sharp}
(y))=\D{\sup_{0\leq t\leq t_0+t_1}} \rho (\varphi^t_{\sharp}
(y))$ for every $x\in A$, hence $f$ is continuous at $x$. It
remains to show the continuity of $f$ at $[L]$. Let $\var$ be a
given positive number.  For every $x\in L$, there is a $t\in \R^+$
such that $\varphi^t (x)\not\in V$. Thus there exists a
neighborhood $U_x$ of $x$ such that $\rho(\varphi^{[0,t]}_{\sharp}
(U_x))<\var$ and $\varphi^t(U_x)\cap V=\varnothing$. If we set
$U=\bigcup_{x\in L} U_x$, then $U$ is a neighborhood of [L] in
$N/L$ with $\rho (\varphi_{\sharp}^{[0,+\infty)} (U))<\var$,
hence $f|_U<\var$. The above argument shows that $f$ is
continuous. Now it is not hand to check that $g:=\int_0^{+\infty}
e^{-t} f(\varphi^t_{\sharp} (x))dt$ is the desired function
\cite{S}. $\square$

\para{Corollary 2.3.} Let $(N,L)$ be a weakly regular index pair for
an isolated invariant set $I$.  Then there exists a subset $L'
\subset N$ such that $L\subset L'$ and $(N,L')$ is a regular
index pair for $I$ .\\ {\it Proof.} Let $g:N/L\lo [0,1]$ be the
Lyapunov function described above. Take $L':= \pi
^{-1}(g^{-1}[0,\var])$ where $\pi :N\lo N/L$ is the natural
projection and $\var \in (0,1)$. Now by Proposition 2.1., $(N,L')$
is a regular index pair for $I$. $\square$

\section{Relative Category}
The relative Lusternik-Schnirelmann category introduced by Fadell
and Husseini \cite{FH} has shown to give important information
about the existence of critical points \cite{Co,My}. For every
topological space $X$ and a closed subset $A\subset X$, the
relative category cat $(X,A)$ is defined to be the minimum $n$
for which there exists an open cover $X=U_0 \cup \ldots \cup U_n$
such that $U_0$ can be retracted to $A$ and $U_i$ is contractible
in $X$ for $1\leq i \leq n$. Unfortunately this definition is not
so efficient in Conley index theory. Indeed in this theory we deal
with a pair in the form of $(N/L,[L])$ where $(N,L)$ is an index
pair for a continuous flow. Now one can see that if $L$ is a
large subset of $N$, then $cat(N/L,[L])$ is possibly a small
number. Since the relative category is used to find a lower bound
for the number of critical points, we can not obtain good
estimates in critical point theory. (See \cite{Po,R1,R2}.) For
this reason, we use a little modification of the relative category
which is more efficient in Conley index theory.

\para{Definition.} (Relative HLS category)
Let $X$ be a topological space and $A\subset X$ be a closed
subset. The relative Homotopy Lusternik-Schnirelmann category
$\nu _H(X,A)$ is defined to be the minimum of $n$ for which there
exists an open covering $X=U_0 \cup\ldots \cup U_n$ such that $A$
is a deformation retract of $U_0$  and $U_i$ is contractible in
$X-A$ for $1\leq i \leq n$.

It is easy to see that $\nu_H(X,A)\geq cat(X,A)$. Unfortunately
the relative category is not easy to compute. It is more
convenient to use cohomological category \cite{MS} to find lower
bounds for it. (See \cite{Co,MO,My,Ru} and references therein for
other tools.) Suppose that $H^*$ is a cohomology functor on $X$. A
subset $S\subset X$ is called cohomologically trivial if the
restriction map $H^k(X)\lo H^k(S)$ is zero for every $k\in \N$.
Similarly for $A\subset S\subset X$, we say that $(S,A)$ is
cohomologically trivial in $(X,A)$ if the restriction map
$H^k(X,A)\lo H^k(S,A)$ is zero for every $k\in \N$.

\para{Definition.} (Relative CLS category) The relative Cohomology
Ljusternik-Schnirelmann category $\nu_C(X,A)$ is defined to be the
minimum $n$ for which there is an open covering $U_0,\cdots ,U_n$
such that $A\subset U_0$ and $(U_0,A)$ is cohomologically trivial
in $(X,A)$ and $U_i \subset X-A$ is cohomologically trivial in
$X-A$ for $1\leq i\leq n$. If such a covering does not exist, we
set $\nu_C(X,A)=+\infty$. When $A=\varnothing$, then
$\nu_C(X):=\nu_C(X,\varnothing)$ is called CLS category.

Recall that the {\em cuplenth} of a topological space $X$ is
defined to be the minimum integer $N>0$ such that for any set of
cohomology classes $\alpha_j\in H^{k_j} (M)$, $j=1,\cdots , N$ of
degree $k_j\geq 1$, the class $(\alpha_1 \cup \cdots \cup
\alpha_N)|_A=j^* (\alpha_1 \cup \cdots \cup \alpha_N)$ is zero.
It is well-known that $\nu_C(X)\geq cuplength(X)$. (cf.
\cite{MS}.) Now we introduce the concept of relative cuplength
which gives a lower bound for the relative CLS category.

\para{Definition.} (Relative Cuplength) For $A\subset X$,
we define the relative cuplength $CL(X,A)$ to be the minimum
integer $N>0$ such that for any set of cohomology classes
$\alpha_j\in H^{k_j} (X)$, $j=1,\cdots , N$ and $\alpha_0\in
H^{k_0} (X,A)$ of degree $k_j>0$, the class $(\alpha_0
\cup\alpha_1 \cup \cdots \cup \alpha_N)$ is zero in $H^*(X,A)$.

\para{Proposition 3.1.} $\nu_C(X,A)\geq CL(X,A)$.\\
{\it Proof.} Let $U$ and $V$ be open subsets of $X$ with $A\subset
V$. Suppose that $\alpha \in H^* (X,A)$ and $\beta \in H^* (X)$
such that $\alpha|_V=0$ in $H^* (V,A)$ and $\beta|_U=0$ in $H^*
(U)$. Since $U$ and $V$ are excisive subsets of $X$, the cup
product map $H^*(X,V)\otimes H^*(X,U)\lo H^* (X,V\bigcup U)$ is
defined \cite{D} and the following diagram is commutative with
exact columns:
$$\begin{array}{ccccc}
H^*(X,V) & \otimes & H^*(X,U) & \lo & H^* (X, V\bigcup U)\\
\downarrow & & \downarrow & & \downarrow\\ H^*(X,A) & \otimes &
H^*(X) & \lo & H^*(X,A)\\ \downarrow & & \downarrow & & \downarrow\\
H^*(V,A) & \otimes & H^*(U) &  & H^*(V\bigcup U,A)
\end{array}$$ A diagram chase shows that $(\alpha \cup \beta)=0$ in $H^*
(V\bigcup U,A)$. Now the proposition is obvious by induction.
$\square$

\para{Example 3.2.}  Let $\pi:E\lo B$ be an $m$-dimensional
orientable vector bundle over a metric space $B$. Then $E$ admits
a Euclidean metric and we may consider the disk bundle $D(E)$ and
sphere bundle $S(E)$. Using the Thom isomorphism, we conclude that
$CL(D(E)/S(E),[S(E)])=CL(D(E),S(E))=cuplength(M)$. (See
\cite{M,MS} for the details.)

\section{Critical Point Theory}

Let $X$ be a complete Finsler manifold and $f\in C^1(X,\R)$. In
\cite{B1,P}, it has been shown that $f$ admits a {\em
pseudo-gradient} vector field i.e. a map $Y:X\lo T(X)$ such that
\begin{itemize}
\item[(i)] The equation $ \dot{x}=Y(x)$ has a unique solution for
every initial point $x_0 \in X$,
\item[(ii)] $Y.f:=\langle Df(x),
Y(x)\rangle \geq \alpha (\| Df(x)\|)$ where $\alpha$ is a
strictly increasing continuous function with $\alpha (0)=0$,
\item[(iii)] $\|Y\|$ is bounded.
\end{itemize}
Therefore we can consider the flow relative to a pseudo-gradient
vector field for $f$. From now on, suppose that $\varphi^t$ is
the flow relative to the pseudo-gradient vector field $Y$ for
$f\in C^1(X,R)$ and $(N,L)$ is an index pair for $\varphi ^t$
such that $f$ is bounded on $V:=\ov{N-L}$. Moreover we assume
that $f$ satisfies Palais-Smale condition in $V$. (See also
\cite{RuS}.)

\para{Definition.} We say that $f$ satisfies Palais-Smale
condition in $V$ if any sequence $\{ x_n \}\subset V$ such that
$f(x_n)$ is bounded and $\|D f(x_n)\|\lo 0$ possesses a convergent
subsequence.

The following lemma shows that under the above assumptions,
$(N,L)$ is a weakly regular index pair. Therefore we can use
Corollary 2.3. to show the existence of a regular index pair in
this situation.

\para{Lemma 4.1.} There exists $T\in\R^+$ such that for every $x\in L$,
 $\varphi^{[0,T]}(x)\not\subset V$. \\ {\it Proof.} Suppose
the contrary, then there is a sequence $x_i\in L$ and $t_i\in
\R^+$ such that $t_i\lo+\infty$ and $\varphi^{[0,t_i]} (x_i)
\subset V$. Since $L$ is positively invariant relative to $N$, we
get $\varphi^{[0,t_i]} (x_i)\subset L$. Now $L\cap V$ is a closed
set which does not contain any critical point of $f$. Since $f$
satisfies Palais-Smale condition in $V$, there is $\delta >0$
such that $Y.f(x)>\delta$ for every $x\in L\cap V$. Therefore
$f(\varphi^{t_i} (x_i))- f(x_i)>\var t_i$. Since $t_i\lo
+\infty$, it follows that $f(\varphi^{t_i} (x_i))\lo +\infty$
which contradicts the boundedness of $f$ on $V$. $\square$

\para{Lemma 4.2.} For large values of $t\in \R^+$,
$(\varphi^t_{\sharp})^{-1} ([L])$ is a neighborhood of $[L]$ in
$N/L$.
\\ {\it Proof.} By the above lemma, there exists $T\in\R^+$
such that $\varphi^{[0,T]} (x)\not\subset V$ for every $x\in L$.
Thus there is a neighborhood $U_x$ of $x$ such that
$\varphi^{[0,T]} (y)\not\subset V$ for every $y\in U_x$. Now
$\bigcup_{x\in L} U_x$ is a neighborhood of $L$, hence it defines
a neighborhood $U$ of $[L]$ with $\varphi^T_{\sharp} (U)=[L]$.
Therefore $(\varphi^t_{\sharp})^{-1} ([L])$ is a neighborhood of
$[L]$ in $N/L$ for every $t\geq T$. $\square$\\

For every $a\in \R$, we define $(N/L)^a=[L]\cup \{x\in N-L|f(x)
\geq a\}.$ The following two lemmas are reformulation of
Deformation Theorems \cite{Ch} for index pairs.

\para{Lemma 4.3.} (First Deformation Theorem) Suppose that $f$ has no
critical points in $(N-L)\cap f^{-1} [a,b]$ for some $a<b$. Then
there exists $T\in\R^+$ such that $(N/L)^a \subset
(\varphi^T_{\sharp})^{-1} ((N/L)^b)$.
\\ {\it Proof.} Since $f$ has no critical points in the closed
subset $V\cap f^{-1} [a,b]$ and $f$ satisfies Palais-Smale
condition in $V$, there is a $\delta>0$ such that $Y.f>\delta$ in
$V\cap f^{-1} [a,b]$. Now for every $T\geq \frac{b-a}{\delta}$ and
$x\in V\cap f^{-1} [a,b]$, we have $f(\varphi^T(x))\geq a+\delta
T\geq b$. $\square$
\para{Lemma 4.4.} (Second Deformation Theorem) Suppose that $f$ has exactly
$m$ critical points $y_1,\ldots, y_m$ in $(N-L)\cap f^{-1}
[c-\var_0, c+\var_0]$ and $f(y_i)=c$ for some $c\in\R$ and $\var
_0 >0$. If $U_i$ is a neighborhood of $x_i$ in $N-L$ and
$T\in\R^+$, then there is an $\var\in[0,\var_0]$ such that
$(N/L)^{c-\var}\subset (\varphi^T_{\sharp})^{-1}
((N/L)^{c+\var})\cup \bigcup_{i=1}^{m} U_i$.  \\ {\it Proof.}
Since $y_i$  is a rest point for $\varphi^t$, there is an open set
$V_i\subset U_i$ such that $\varphi^{[0,T]} (x) \cap
V_i=\varnothing$ for every $x\not\in U_i$. Since $f$ has no
critical point in  $V\cap f^{-1}
[c-\var_0,c+\var_0]-\bigcup_{i=1}^n V_i$ which is a closed subset
of $V$ and $f$ satisfies Palais-Smale condition in $V$, there
exists $\delta>0$ such that $Y.f>\delta$ in $V\cap f^{-1}
[c-\var_0,c+\var_0]-\bigcup_{i=1}^n V_i$. If we set
$\var=\frac{T\delta}{2}$, then for every $x\in V\cap f^{-1}
[c-\var, c+\var]-\bigcup_{i=1}^n U_i$, we have
$f(\varphi^T(x))\geq f(x)+T\delta\geq c-\var+2\var=c+\var$.
$\square$
\para{Proposition 4.5.} Suppose that $f$ has a finite number of critical
points $x_1,\ldots , x_n$ in $V$ and $U_i$ is a neighborhood of
$x_i$ in $N-L$ for $1\leq i \leq n$. Then there are
$t_0,t_1,\ldots,t_n\in \R^+$ such that
$N/L=(\varphi^{t_0}_{\sharp})^{-1} ([L])\cup \bigcup_{i=1}^n
(\varphi^{t_i}_{\sharp})^{-1} (U_i)$.  In particular if $(N,L)$ is
a regular index pair, then $N=(\varphi^{t_0}_{\natural})^{-1}
(L)\cup \bigcup_{i=1}^n (\varphi^{t_i}_{\natural})^{-1} (U_i)$\\
{\it Proof.} We use induction on the number of critical values of
$f|_V$. For $k=0$, $f$ has no critical points in $V$. Since $f$ is
bounded on $V$, there are $a<b$ such that $(N/L)^b=[L]$ and
$(N/L)^a=N/L$. Thus by Lemma 4.3.,
$(\varphi_{\sharp}^{T})^{-1}([L])=N/L$ for some $T\in \R^+$. Now
suppose that $f$ has $k+1$ critical values $c_0<c_1<c\ldots <c_k$
and $x_{\ell}, \ldots , x_n$ are critical points with $f(x_i)=c_k$
for $\ell \leq i \leq n$. If we use Lemma 4.4. for
$\var_0=\frac{c_k-c_{k-1}}{2}$ and $T=1$, we obtain an
$\var\in[0,\var _0]$ such that $(N/L)^{c-\var}\subset \varphi
_{\sharp}^{-1} ((N/L)^{c+\var}) \cup \bigcup_{i=\ell}^n U_i$.
Moreover by Lemma 4.3., there is a $T\in \R^+$ such that
$(N/L)^{c-\var}\subset \varphi _{\sharp} ^{-T}([L])$. Therefore we
have $(N/L)^{c-\var}\subset (\varphi^{T+1}_{\sharp})^{-1}
([L])\cup \bigcup_{i=\ell}^{n} U_i$. Now if we set $L_1=N\cap
f^{-1} [c-\var,+\infty)$, then $(N,L_1)$ is an index pair with $k$
critical values, hence there are $t'_0,t'_1,\ldots, t'_{\ell-1}$
such that $N/L_1=(\varphi_{\sharp}^{t'_0})^{-1} ([L_1])\cup
\bigcup_{i=1}^{\ell -1} (\varphi^{t'_i}_{\sharp})^{-1} (U_i)$. Now
we set $t_0=t'_0+T+1, t_i=t'_i$ for $1\leq i \leq \ell -1$ and
$t_i=t'_0$ for $\ell \leq i\leq n$. Then it is not hard to check
that $$N/L=(\varphi_{\sharp}^{t_0})^{-1} ([L])\cup \bigcup_{i=1}^n
(\varphi^{t_i}_{\sharp})^{-1} (U_i).\  \square $$

\para{Theorem 4.6.} $f$ has at least $\nu_H(N/L,[L])$ critical points
in $V$. Moreover if $(N,L)$ is a regular index pair, then $f$ has
at least $\nu_H(N,L)$ critical points in $V$.\\ {\it Proof.} We
may assume that $f$ has a finite number of critical points
$x_1,\ldots, x_n$ in $V$. Since $X$ is a Banach manifold, we may
choose a contractible open set $U_i\subset N-L$ which contains
$x_i$ for $1\leq i\leq n$. Now by the above proposition, there are
$t_0,t_1,\ldots, t_n\in\R^+$ such that
$N/L=(\varphi_{\sharp}^{t_0})^{-1} ([L])\cup \bigcup_{i=1}^n
(\varphi^{t_i}_{\sharp})^{-1} (U_i)$. It is easy to see that each
$(\varphi^{t_i}_{\sharp})^{-1} (U_i)$ is contractible in $N-L$.
The only problem is that $(\varphi^{t_0}_{\sharp})^{-1} ([L])$ is
not an open set. By Lemma 4.2, there is $T\in \R^+$ such that
$int ((\varphi^T_{\sharp})^{-1}([L]))\neq \varnothing$. Now we
replace $(\varphi^{t_0}_{\sharp})^{-1}([L])$ by
$int((\varphi^{t_0+T}_{\sharp})^{-1} ([L]))$ which is obviously
invariant under $\varphi_{\sharp}^t$ and conrtactible to $[L]$.
Therefore we obtain an open covering which follows that
$\nu_H(N/L,[L])\leq n$. A similar argument for $\varphi
_{\natural}$ shows that if $(N,L)$ is a regular index
 pair, then $\nu_H(N,L)\leq n$.  $\square$

\para{Example 4.7.} Let $E$ be an orientable vector bundle over a
Finsler manifold $X$ and $B$ be a closed subset of $X$. We denote
the disk and sphere bundles of $E|_B$ by $D$ and $S$
respectively. Let $f:E\lo \R$ be a $C^1$ function satisfying
Palai-Smale condition in $B$ and $\varphi^t$ be the flow relative
to a pseudo-gradient vector field for $f$ such that $(D,S)$ is an
index pair for $\varphi^t$. Then by Example 2.3. and the above
theorem, $f$ has at least $\nu_H(D,S)\geq cuplength(B)$ critical
points in $D$. If we consider the case $X=M\times \R^n,\
E=X\times \R^m$ and $B=M\times D^n\subset M\times \R^n$, then we
obtain a noncompact version of a well-known result of \cite{CZ1}.

\para{Acknowledgment.} The author would like to thank Institute for
Studies in Theoretical Physics and Mathematics for supporting this
research.

\end{document}